\documentclass[psamsfonts]{amsart}  

\usepackage[T1]{fontenc}
\usepackage{amssymb}
\usepackage{amsthm}
\usepackage{amsmath}

\usepackage[dvips]{graphicx}

\usepackage{graphicx}
\usepackage{amssymb}                       
\usepackage{amsmath}
\usepackage{color}
\usepackage{times}

\usepackage[unicode,bookmarks,colorlinks]{hyperref}
\hypersetup{
   linkcolor=brickred,
}

\definecolor{mahogany}{cmyk}{0, 0.77, 0.87, 0}
\definecolor{salmon}{cmyk}{0, 0.53, 0.38, 0}
\definecolor{melon}{cmyk}{0, 0.46, 0.50, 0}
\definecolor{yellowgreen}{cmyk}{0.44, 0, 0.74, 0}
\definecolor{brickred}{cmyk}{0, 0.89, 0.94, 0.28}
\definecolor{OliveGreen}{cmyk}{0.64, 0, 0.95, 0.40}
\definecolor{RawSienna}{cmyk}{0, 0.72, 1.0, 0.45}
\definecolor{ZurichRed}{rgb}{1, 0, 0} 

\usepackage{fancyhdr}
\usepackage{amsmath,amstext,amssymb,amsopn,amsthm}
\usepackage{amsmath,amssymb,amsthm}
\usepackage[mathscr]{eucal}

\DeclareMathOperator{\sign}{sign}
\newcommand{\tr}{\operatorname{tr}}

\reversemarginpar

\definecolor{rb}{rgb}{0.1,0.2, 0.7}

\begin{document}

\newtheorem{lemma}[thm]{Lemma}
\newtheorem{proposition}{Proposition}
\newtheorem{theorem}{Theorem}[section]
\newtheorem{deff}[thm]{Definition}
\newtheorem{case}[thm]{Case}
\newtheorem{prop}[thm]{Proposition}
\newtheorem{example}{Example}

\newtheorem{corollary}{Corollary}

\theoremstyle{definition}
\newtheorem{remark}{Remark}

\numberwithin{equation}{section}
\numberwithin{definition}{section}
\numberwithin{corollary}{section}

\numberwithin{theorem}{section}

\numberwithin{remark}{section}
\numberwithin{example}{section}
\numberwithin{proposition}{section}

\newcommand{\gap}{\lambda_{2,D}^V-\lambda_{1,D}^V}
\newcommand{\gapR}{\lambda_{2,R}-\lambda_{1,R}}
\newcommand{\bD}{\mathrm{I\! D\!}}
\newcommand{\calD}{\mathcal{D}}
\newcommand{\calA}{\mathcal{A}}

\newcommand{\conjugate}[1]{\overline{#1}}
\newcommand{\abs}[1]{\left| #1 \right|}
\newcommand{\cl}[1]{\overline{#1}}
\newcommand{\expr}[1]{\left( #1 \right)}
\newcommand{\set}[1]{\left\{ #1 \right\}}
\newcommand{\Rp}{\mathbb R^{d+1}_+}

\newcommand{\calC}{\mathcal{C}}
\newcommand{\calE}{\mathcal{E}}
\newcommand{\calF}{\mathcal{F}}
\newcommand{\Rd}{\mathbb{R}^d}
\newcommand{\BR}{\mathcal{B}(\Rd)}
\newcommand{\R}{\mathbb{R}}
\newcommand{\T}{\mathbb{T}}
\newcommand{\D}{\mathbb{D}}

\newcommand{\al}{\alpha}
\newcommand{\RR}[1]{\mathbb{#1}}
\newcommand{\bR}{\mathrm{I\! R\!}}
\newcommand{\ga}{\gamma}
\newcommand{\om}{\omega}
\newcommand{\A}{\mathbb{A}}
\newcommand{\bH}{\mathbb{H}}

\newcommand{\bb}[1]{\mathbb{#1}}
\newcommand{\bI}{\bb{I}}
\newcommand{\bN}{\bb{N}}

\newcommand{\uS}{\mathbb{S}}
\newcommand{\M}{{\mathcal{M}}}
\newcommand{\calB}{{\mathcal{B}}}

\newcommand{\W}{{\mathcal{W}}}

\newcommand{\m}{{\mathcal{m}}}

\newcommand {\mac}[1] { \mathbb{#1} }

\newcommand{\bC}{\Bbb C}

\newcommand\mb{\mbox{$\frac{\pi}{2}-\theta$}}
\newcommand\pie{\frac{\pi}{2}-\frac{\pi}{2p}}

\newtheorem{rem}[theorem]{Remark}
\newtheorem{dfn}[theorem]{Definition}
\theoremstyle{definition}
\newtheorem{ex}[theorem]{Example}
\numberwithin{equation}{section}

\newcommand{\Pro}{\mathbb{P}}
\newcommand\F{\mathcal{F}}
\newcommand\E{\mathbb{E}}
\newcommand\e{\varepsilon}
\def\H{\mathcal{H}}
\def\t{\tau}

\title[Martingale inequalities]{A dual approach to $L^p$ estimates for orthogonal martingales satisfying the differential  subordination}

\author{Rodrigo Ba\~nuelos}
\address{Department of Mathematics, Purdue University, West Lafayette, IN 47907, USA}
\email{banuelos@math.purdue.edu}

\author{Adam Os\k{e}kowski}
\thanks{A. Os\k{e}kowski is supported by Narodowe Centrum Nauki (Poland), grant no. DEC-2014/14/E/ST1/00532.}\address{Faculty of Mathematics, Informatics and Mechanics\\
 University of Warsaw\\
Banacha 2, 02-097 Warsaw\\
Poland}
\email{A.Osekowski@mimuw.edu.pl}


\begin{abstract}
We develop a new dual approach, which can be used to establish sharp $L^p$ bounds for orthogonal differentially subordinate martingales. The method rests on the construction of a certain special function enjoying appropriate majorization properties and  a novel differential condition.
\end{abstract}

\maketitle

\section{Introduction}

Sharp martingale and semimartingale inequalities play an important role in applications in functional and harmonic analysis, as they often lead to  tight bounds for wide classes of operators, including maximal functions, Fourier multipliers and sparse forms. There is a powerful technique, called the Bellman function method, or Burkholder's method, which can be used to study such results for quite general stochastic processes. The idea behind the approach can be formulated as follows: the validity of a given estimate is equivalent to the existence of a certain special function, enjoying appropriate majorization and concavity properties. It should be emphasized that the method brings an additional insight into the estimate under investigation: not only does the special function yield the bound, but it also encodes the structure of extremizers (i.e., processes for which equality, or almost equality, is attained) and often embeds the result into a wider class of inequalities. Furthermore, the special functions are of interest in their own right, as they might exhibit deep connections between several seemingly unrelated contexts such as the interplay between Burkholder-type functionals, nonlinear Cauchy-Riemann equations and the theory of quasiconvex functions, see for example \cite{AIP,BMS,I}.

From the historical perspective, the Bellman function technique originates in the theory of optimal control and the foundational works of Bellman from the fifties \cite{Be}. Its connection to the theory of martingales appeared  for the first time in the eighties, in the seminal paper \cite{B0} of Burkholder, who used the method to obtain a wide class of sharp inequalities for martingale transforms (and applied them in the study of the Haar system). This path was further explored by Burkholder and his students in the semimartingale case (see the monograph \cite{Os} for the overview of the results). There is also the analytic branch of the approach, developed in the mid-nineties by Nazarov, Treil and Volberg \cite{NT,NTV}, who showed that the method can be exploited beyond the probability theory. Since then, the technique has been exploited in numerous works on various topics, such as maximal estimates, fractional estimates, bounds for $BMO$ functions and martingales, weighted theory, and sparse domination which connects further with the Calder\'on-Zygmund theory.

In general, the Bellman function leading to a given estimate under investigation need not be unique. In fact,  a clever choice often enables to avoid, or at least reduce the number of complicated calculations, which arise during the verification of properties of the special functions. However, this non-uniqueness phenomenon goes much deeper. Namely, some estimates can be proven with the use of completely different Bellman setups (i.e., with the use of functions for which the concavity-type requirement is entirely different), each of which has its advantages and deficiencies. As an illustration, let us present a specific example: let us briefly discuss $L^p$ inequalities for differentially subordinate martingales.

\smallskip

We start with some background and notation. Suppose that $(\Omega,\F,\mathbb{P})$ is a probability space equipped with a discrete-time filtration. Let $f=(f_n)_{n\geq 0}$, $g=(g_n)_{n\geq 0}$ be two adapted martingales, with difference sequences $df=(df_n)_{n\geq 0}$, $dg=(dg_n)_{n\geq 0}$ given by
$$ df_0=f_0,\quad df_n=f_n-f_{n-1},\qquad n=1,\,2,\,\ldots,$$
and similarly for $dg$. We say that $g$ is differentially subordinate to $f$, if we have
\begin{equation}\label{differ} 
|dg_n|\leq |df_n| \qquad \mbox{almost surely for }n=0,\,1,\,2,\,\ldots.
\end{equation}
The following  is the celebrated $L^p$ estimate of Burkholder \cite{B0}.  We use the notation $p^*=\max\{p,p'\}$, where $p'=p/(p-1)$ is the harmonic conjugate to $p$, and denote $\|f\|_p=\sup_n \|f_n\|_p$.

\begin{theorem}\label{B_main}
If $f$, $g$ are martingales such that $g$ is differentially subordinate to $f$, then
\begin{equation}\label{moment}
 \|g\|_p\leq (p^*-1)\|f\|_p,\qquad 1<p<\infty.
\end{equation}
For each $p$ the constant $p^*-1$ is the best possible: it cannot be replaced by a smaller number.
\end{theorem}

It is a very important result for applications: it leads to tight estimates for a wide class of Fourier multipliers (including second-order Riesz transforms and the Beurling-Ahlfors operator \cite{BMH,BW,BW2,NV}). 
Probably the simplest proof of this statement (which will be referred to as ``the direct proof'') comes from Burkholder's survey \cite{B1} and proceeds as follows. Suppose that the function $U:\R^2\to \R$ satisfies 

\medskip

(i)  $U(x,y)\leq 0$ if $|y|\leq |x|$,

\smallskip

(ii) $U(x,y)\geq |y|^p-(p^*-1)^p|x|^p$ for all $x,\,y\in \R$,

\smallskip

(iii) $U$ is concave along any line contained in $\R^2$, of slope belonging to $[-1,1]$.

\medskip

\noindent Then for any pair $f$, $g$ as in the statement of Theorem \ref{B_main}, the property (iii) implies that the composition $(U(f_n,g_n))_{n\geq 0}$ is a supermartingale.  Combining this with the other two conditions, one obtains, for any $n$,
$$ \E\Big[|g_n|^p-(p^*-1)^p|f_n|^p\Big]\leq \E U(f_n,g_n)\leq \E U(f_0,g_0)\leq 0.$$
This yields \eqref{moment}, by letting $n\to \infty$. Burkholder showed that the function
$$
 U(x,y)=p\left(1-\frac{1}{p^*}\right)^{p-1}(|y|-(p^*-1)|x|)(|x|+|y|)^{p-1}
$$
has all the required properties. It is worth mentioning that the choice of the special function is not unique. For instance, for $p>2$, the function
$$ \tilde{U}(x,y)=\begin{cases}
U(x,y) & \mbox{if }|y|>(p-1)|x|,\\
|y|^p-(p-1)^p|x|^p & \mbox{if }|y|\leq (p-1)|x|
\end{cases}$$
also does the job (and so does any convex combination of $U$ and $\tilde{U}$). 

\smallskip

There is a dual method of proving \eqref{moment}, which was developed and simplified in several works by Nazarov, Treil and Volberg \cite{NT,NTV}. The argument also exploits a certain special function: suppose that a $C^2$ function $B:\R^2\to \R$ satisfies

\medskip

(i) $B(x,z)\geq |xz|$ for all $x,\,z\in \R$.

\smallskip

(ii) $B(x,z)\leq \frac{K^p|x|^p}{p}+\frac{L^{p'}|z|^{p'}}{p'}$ for all $x,\,z\in \R$.

\smallskip

(iii) For any $x,\,z,\,h,\,k\in \R$ we have
\begin{equation}\label{conv0}
B_{xx}(x,z)h^2+2B_{xz}(x,z)hk+B_{zz}(x,z)k^2\geq 2|h||k|.
\end{equation}

\smallskip

\noindent Then \eqref{moment} holds with the constant $KL$. To see this, pick an arbitrary pair $(f,g)$ as in the statement and an arbitrary $L^{p'}$-bounded, a priori unrelated martingale $h$. It can be shown that (iii) implies the inequality
$$ \E \Big[B(f_n,h_n)-B(f_{n-1},h_{n-1})\Big]\geq \E |df_n||dh_n|,\qquad n=1,\,2,\,\ldots.$$
Hence, using the other two conditions, we obtain
\begin{equation}\label{chain1}
\begin{split}
 \E \sum_{k=0}^n |df_k||dh_k|&\leq \E \left[B(f_0,h_0)+\sum_{k=1}^n |df_k||dh_k|\right]\\
&\leq \E B(f_n,h_n)\leq \E\left(\frac{K^p|f_n|^p}{p}+\frac{L^{p'}|h_n|^{p'}}{p'}\right).
\end{split}
\end{equation}
However, by the differential subordination and the orthogonality of martingale differences,
\begin{equation}\label{chain2}
 \E \sum_{k=0}^n |df_k||dh_k|\geq \E \sum_{k=0}^n |dg_k||dh_k|\geq \E \sum_{k=0}^n dg_kdh_k=\E g_nh_n.
\end{equation}
Plugging this above and using duality (recall that $h$ was arbitrary), we obtain $\|g_n\|_p\leq KL \|f_n\|_p$ and it remains to let $n\to \infty$ to get \eqref{moment} (with $p^*-1$ replaced by $KL$). Nazarov and Treil \cite{NT} constructed the function
$$ B(x,z)=|x|^p+|z|^{p'}+\begin{cases}
|x|^2|z|^{2-p'} & \mbox{if }|x|^p\leq |z|^{p'},\\
\displaystyle \frac{2}{p}|x|^p+\left(\frac{2}{p'}-1\right)|z|^{p'} & \mbox{if }|x|^p\geq |z|^{p'},
\end{cases}$$
which enabled them to prove \eqref{moment} with a suboptimal constant. 
 The special function leading to the sharp $L^p$ bound was invented by Ba\~nuelos, Ga\l \k azka and Os\k ekowski \cite{BGO}: the description is a bit technical and will not be presented here. Let us conclude this discussion by the short comparison of both approaches. The direct proof of Burkholder is less technical and more general: it can be easily adjusted to the study of related inequalities (cf. \cite{Os}). On the other hand, the dual method has the advantage of leading directly to many interesting applications in harmonic analysis and semigroup theory: see e.g. \cite{DC,DV,Wr}.

\smallskip

The purpose of this paper is to develop a dual approach to a different class of inequalities. We will be interested in the sharp $L^p$ bounds for differentially subordinate martingales satisfying the additional orthogonality assumption. We need to change the setup to the continuous-time case, since the orthogonality condition forces the dominated martingale to have continuous paths  (see Lemma \ref{l1W} below). From now on, we assume that $(\Omega,\F,\mathbb{P})$ is a complete probability space, equipped with a filtration $(\F_t)_{t\geq 0}$ such that $\F_0$ contains all the events of probability $0$. Let $X=(X_t)_{t\geq 0}$, $Y=(Y_t)_{t\geq 0}$ be two real-valued local martingales, which have right-continuous paths with limits from the left. That is, for almost all $\omega\in \Omega$, the functions $t\mapsto X_t(\omega)$ and $t\mapsto Y_t(\omega)$ are right-continuous on $[0,\infty)$ and have left-sided limit at each $t\in (0,\infty)$.  The local martingales $X$, $Y$ are said to be orthogonal, if $[X,Y]$ is a constant process. Furthermore, $Y$ is differentially subordinate to $X$, if the difference $([X,X]_t-[Y,Y]_t)$ is nonnegative and nondecreasing. Here and below, the symbol $[X,Y]=([X,Y]_t)_{t\geq 0}$ will stand for the quadratic variation between $X$ and $Y$,  defined as the limit in probability
\begin{equation}\label{[X,Y]}
 [X,Y]_t=\operatorname*{\limsup}_{n\to \infty} \sum_{j=1}^{k_n} (X_{t_j^{(n)}}-X_{t_{j-1}^{(n)}})(Y_{t_j^{(n)}}-Y_{t_{j-1}^{(n)}}),
\end{equation}
where $(t_j^{(n)})_{j=0}^{k_n}$ is an arbitrary sequence of partitions $0=t_0^{(n)}<t_1^{(n)}<t_2^{(n)}<\ldots<t_{k_n}^{(n)}=t$ of $[0,t]$ with a diameter tending to $0$. 
It is not difficult to see that the notion of continuous-time differential subordination is a generalization of \eqref{differ}.

Recall the following result of Ba\~nuelos and Wang from \cite{BW}. In analogy with the discrete-time setting, we use the notation $\|X\|_p=\sup_{t\geq 0}\|X_t\|_p$ for all $p$.

\begin{theorem}
Suppose that $X$, $Y$ are orthogonal local martingales such that $Y$ is differentially subordinate to $X$. Then we have the sharp estimate
\begin{equation}\label{mainin}
\|Y\|_p\leq \cot\left(\frac{\pi}{2p^*}\right)\|X\|_p,\qquad 1<p<\infty.
\end{equation}
\end{theorem}

This result has importance consequences in harmonic analysis, as it leads to sharp $L^p$ bounds for the first-order Riesz transforms, in the Euclidean and non-Euclidean settings \cite{BO,BW}. As for the proof, it turns out that the estimate \eqref{mainin} follows from the existence of a function $U:\R^2\to \R$ satisfying 

\smallskip

(i) $U(x,y)\leq 0$ for $|y|\leq |x|$.

\smallskip

(ii) $U(x,y)\geq |y|^p-\big(\cot(\pi/(2p^*))\big)^p|x|^p$ for all $x,\,y\in \R$.

\smallskip

(iii) $U$ is superharmonic.

\smallskip

As proved in \cite{BW}, the explicit formula for $U$ is conveniently expressed in terms of polar coordinates $|x|=R\cos \theta$, $|y|=R\sin \theta$, $\theta\in [0,\pi/2]$. Namely, for $1<p\leq 2$ one lets
$$ U(x,y)=-\frac{\sin^{p-1}(\pi/2p)}{\cos(\pi/2p)}R^p\cos(p\theta),$$
while for $p>2$ we can take
$$ U(x,y)=\begin{cases}
\displaystyle \frac{\cos^{p-1}(\pi/2p)}{\sin(\pi/2p)}R^p\cos(p(\mb)) & \mbox{if }\theta>\frac{\pi}{2}-\frac{\pi}{2p},\\
|y|^p-(\cot\frac{\pi}{2p})^p|x|^p & \mbox{if } \theta\leq \frac{\pi}{2}-\frac{\pi}{2p}.
\end{cases}$$

The above proof can  be regarded as a direct argument. This gives rise to the very interesting and natural question of whether the estimate \eqref{mainin} can be established using  an alternative dual approach. One of the contributions of this work is to give the affirmative answer to this question. As we will see, the argument will have some common features with that described above in the context of discrete-time differential subordination, but it will also include an extra technical step which enables the handling of the orthogonality condition. From now on, we assume that $p\neq 2$: the case $p=2$ is already covered by the preceding context of differential subordination (i.e., without the orthogonality assumption).

\smallskip
Let us say a few words about the structure of the remaining part of the paper. The next section contains the description of the dual approach to \eqref{mainin}. In Section 3 we exploit the method to establish the $L^p$ bound. The final part of the paper is devoted to the following application: we show how the dual Bellman function leads directly, with no use of martingale methods, to sharp $L^p$ bounds for Riesz transforms on $\R^d$.

\section{On the approach}

We start with a structural property of local martingales satisfying the differential subordination and orthogonality. First, recall the following well-known fact in stochastic analysis (see e.g. \cite[Chapters VI and VII]{DM} for details). Namely, for any local martingale $X$ there exists a unique continuous local martingale part $X^c$ of $X$ satisfying
$$ [X,X]_t=|X_0|^2+[X^c,X^c]_t+\sum_{0<s\leq t} |\Delta X_s|^2$$
for all $t\geq 0$. Here $\Delta X_s=X_s-X_{s-}$ denotes the jump of $X$ at time $s$. Furthermore, we have that  $[X^c,X^c]=[X,X]^c$, the pathwise continuous part of $[X,X]$. Here is Lemma 2.1 from \cite{BW2}.

\begin{lemma}\label{l1W}
If $X$ and $Y$ are local martingales, then $Y$ is differentially subordinate and orthogonal to $X$ if and only if $Y^c$ is differentially subordinate and orthogonal to $X^c$, $|Y_0|\leq |X_0|$ and $Y$ has continuous paths.
\end{lemma}

We are ready to describe our approach to $L^p$ estimates for orthogonal differentially subordinate local martingales. In what follows, $D^2B$ denotes the Hessian matrix of $B$.

\begin{proposition}\label{main_tool}
Let $1<p<\infty$, $b>0$, $\alpha>0$ and $\beta>0$ be fixed parameters.  
Let $B:\R^2\to \R$ be a convex $C^1$ function. We assume further that $B$ is of class $C^2$ on $D_i$, $i\geq 1$, where $(D_i)_i$ is a sequence of open connected sets such that $\bigcup_i \overline{D_i}=\R^2$. If $B$ satisfies the requirements
\begin{equation}\label{initial}
B(x,z)\geq b|xz| \qquad \mbox{for all }(x,z)\in \R^2,
\end{equation}
\begin{equation}\label{majorization}
B(x,z)\leq \alpha |x|^p+\beta|z|^{p'}\qquad \mbox{for all }(x,z)\in \R^2
\end{equation}
and
\begin{equation}\label{det2}
\operatorname*{det}D^2B(x,z)\geq b^2 \qquad \mbox{for all }(x,z)\in \bigcup_i D_i,
\end{equation}
then for any orthogonal local martingales $X$, $Y$ such that $Y$ is differentially subordinate to $X$ we have
\begin{equation}\label{Lp}
 \|Y\|_p\leq \frac{(p\alpha)^{1/p}(p'\beta)^{1/p'}}{b}\|X\|_p
\end{equation}
and
\begin{equation}\label{Lp'}
 \|Y\|_{p'}\leq \frac{(p\alpha)^{1/p}(p'\beta)^{1/p'}}{b}\|X\|_{p'}.
\end{equation}
\end{proposition}

\begin{remark}\label{dual_remark}
Here we encounter an unusual phenomenon. The existence of a function $B$ satisfying \eqref{initial}, \eqref{majorization} and \eqref{det2} yields the validity of both estimates in $L^p$ and $L^{p'}$. This will allow us to restrict ourselves, while showing \eqref{mainin}, to the case $p>2$. In contrast, the direct proof of \eqref{moment} or \eqref{mainin} requires the construction of special functions for each exponent.
\end{remark}

\begin{proof}[Proof of Proposition \ref{main_tool}]
The idea is to apply It\^o's formula to the process $B(X,Z)$ and to analyze carefully the appearing terms. However, since $B$ does not need to be of class $C^2$ on its whole domain, we will 
 use a mollification argument with allows to approximate $B$ with suitable smooth functions. The proof is quite long, and we have decided to split it into several intermediate steps.

\smallskip

\emph{Step 1. Some properties of $B$, a mollification.} 
The convexity of $B$, together with the condition \eqref{det2}, imply that $B_{zz}> 0$ on  $\bigcup_i D_i$ (this will guarantee that certain expressions below will make sense: the term $B_{zz}$ will appear in several denominators). Fix $(x,z)\in \bigcup_i D_i$, let $a(x,z)=B_{xz}(x,z)/B_{zz}(x,z)$ and observe that we have the estimate
\begin{equation}\label{conv}
\begin{split}
 &B_{xx}(x,z)h^2+2B_{xz}(x,z)hk+B_{zz}(x,z)k^2\\
&\qquad \qquad  \geq \frac{b^2}{B_{zz}(x,z)}\cdot h^2+B_{zz}(x,z)(k+a(x,z)h)^2
 \end{split}
\end{equation}
 for all $h,\,k\in \R$. Indeed, since $B_{zz}>0$, this is equivalent to $(\operatorname*{det}D^2B(x,z)-b^2)h^2\geq 0$. 

\smallskip

 Next, let $g:\R^2\to[0,\infty)$ be a $C^\infty$ function, supported on a unit ball of $\R^2$ and satisfying $\int_{\R^2} g=1.$ For a given positive integer $\ell$, define $B^\ell:\R^2\to \R$ by the convolution
\begin{equation}\label{Mollified:1} B^\ell(x,z)=\int_{\R^2} B(x+ u/\ell,z+v/\ell)g(u,v)\mbox{d}u\mbox{d}v.
\end{equation}
This function is of class $C^\infty$ and inherits the convexity. Using integration by parts, we get
\begin{equation}\label{i1}
 B^\ell_{xx}(x,z)=\int_{\R^2} B_{xx}(x+u/\ell,z+v/\ell)g(u,v)\mbox{d}u\mbox{d}v,
 \end{equation}
 \begin{equation}\label{i2} B^\ell_{xz}(x,z)=\int_{\R^2} B_{xz}(x+u/\ell,z+v/\ell)g(u,v)\mbox{d}u\mbox{d}v,
 \end{equation}
\begin{equation}\label{i3} B^\ell_{zz}(x,z)=\int_{\R^2} B_{zz}(x+u/\ell,z+v/\ell)g(u,v)\mbox{d}u\mbox{d}v.
\end{equation}
Now we apply \eqref{conv} to the point $(x+u/\ell,z+v/\ell)\in \bigcup_i D_i$. Multiplying both sides by $g(u,v)$ and integrating over $\R^2$ with respect to $(u,v)$, we obtain
\begin{align*} &B_{xx}^\ell(x,z)h^2+2B_{xz}^\ell(x,z)hk+B_{zz}^\ell(x,z)k^2\\
&\qquad \geq \int_{\R^2}\bigg[\frac{b^2}{B_{zz}(x+u/\ell,z+v/\ell)}\cdot h^2\\
&\qquad \qquad \qquad \qquad +B_{zz}(x+u/\ell,z+v/\ell)\big(k+a(x+u/\ell,z+v/\ell)h\big)^2\bigg] g(u,v)\mbox{d}u\mbox{d}v,
\end{align*}
by the identities \eqref{i1}, \eqref{i2} and \eqref{i3}. 
\smallskip

\emph{Step 2. The application of the It\^o formula.} 
Now, fix a pair $(X,Y)$ of orthogonal local martingales such that $Y$ is differentially subordinate to $X$. We may assume that $X$ is bounded in $L^p$, since otherwise the claim is obvious. Let $Z$ be another adapted local martingale (a priori, no relations to $X$, $Y$ are assumed), satisfying $\|Z\|_{p'}<\infty$. Let $(\tau_n)_{n\geq 0}$ be a sequence of stopping times, given by
\begin{equation}\label{def_stop}
 \tau_n=\inf\big\{t\geq 0\,:\,|X_t|+|Z_t|\geq n\big\}.
\end{equation}
Since $B^\ell$ is of class $C^\infty$, we may apply the It\^o formula to the process $(B^\ell(X_t,Z_t))_{t\geq 0}$, obtaining
\begin{equation}\label{ito}
 I_0=I_1+I_2+\frac{1}{2}I_3+I_4,
\end{equation}
where
\begin{align*}
I_0&=B^\ell(X_{\tau_n\wedge t},Z_{\tau_n\wedge t}),\\
I_1&=B^\ell(X_0,Z_0),\\
I_2&=\int_0^{\tau_n\wedge t} B_x^\ell(X_{s-},Z_{s-})\mbox{d}X_s+\int_0^{\tau_n\wedge t} B_z^\ell(X_{s-},Z_{s-})\mbox{d}Z_s,\\
I_3&=\sum_{s\leq \tau_n\wedge t} \!\!\Big(B^\ell(X_s,Z_s)-B^\ell(X_{s-},Z_{s-})-B_x^\ell(X_{s-},Z_{s-})\Delta X_s-B_z^\ell(X_{s-},Z_{s-})\Delta Z_s\Big),\\
I_4&=\int_0^{\tau_n\wedge t} B_{xx}^\ell(X_{s-},Z_{s-})\mbox{d}[X^c,X^c]_s+2\int_0^{\tau_n\wedge t} B_{xz}^\ell(X_{s-},Z_{s-})\mbox{d}[X^c,Z^c]_s\\
&\quad +\int_0^{\tau_n\wedge t} B_{zz}^\ell(X_{s-},Z_{s-})\mbox{d}[Z^c,Z^c]_s.
\end{align*}

\smallskip

\emph{Step 3. The analysis of the terms $I_0$-$I_3$}. First, note that the majorization \eqref{majorization} implies 
$$ I_0\leq \int_{\R^2} \Big(\alpha|X_{\tau_n\wedge t}+u/\ell|^p+\beta|Z_{\tau_n\wedge t}+v/\ell|^{p'}\Big)g(u,v)\mbox{d}u\mbox{d}v.$$ 
We will handle the term $I_1$ in the next step, when we carry over the limiting procedure $\ell\to \infty$. By the definition of $\tau_n$, both stochastic integrals in $I_2$ are martingales and hence have vanishing expectation. The term $I_3$ is positive, since $B^\ell$ is convex. To deal with $I_4$, fix two stopping times $\sigma\leq \eta$. Using the notation $[S,T]_s^u=[S,T]_u-[S,T]_s$ and approximating square brackets by sums as in \eqref{[X,Y]}, we see that the last inequality of Step 1 implies
\begin{equation*}
\begin{split}
& B_{xx}^\ell(X_{\sigma-},Z_{\sigma-})[X^c,X^c]_{\sigma}^{\eta}
+2B_{xz}^\ell(X_{\sigma-},Z_{\sigma-})[X^{c},Z^{c}]_{\sigma}^{\eta}+
B_{zz}^\ell(X_{\sigma-},Z_{\sigma-})[Z^{c},Z^{c}]_{\sigma}^{\eta}\\
& \geq \int_{\R^2}\bigg[\frac{b^2}{B_{zz}(X_{\sigma-}+u/\ell,Z_{\sigma-}+v/\ell)}\cdot [X^c,X^c]_\sigma^\eta\\
&\quad \qquad +B_{zz}(X_{\sigma-}+u/\ell,Z_{\sigma-}+v/\ell)\big[Z^c+(\xi^{(u,v)})^c,Z^c+(\xi^{(u,v)})^c\big]_\sigma^\eta\bigg] g(u,v)\mbox{d}u\mbox{d}v.
\end{split}
\end{equation*}
Here $\xi^{(u,v)}$ denotes the local martingale
$$ \xi^{(u,v)}_t=a(X_{\sigma-}+u/\ell,Z_{\sigma-}+v/\ell)(X_{\eta\wedge t}^c-X_{\sigma\wedge t}^c).$$
By Lemma \ref{l1W},  we have $[X^c,X^c]_{\sigma}^\eta\geq [Y^c,Y^c]_{\sigma}^\eta$, so
\begin{equation*}
\begin{split}
& B_{xx}^\ell(X_{\sigma-},Z_{\sigma-})[X^c,X^c]_{\sigma}^{\eta}
+2B_{xz}^\ell(X_{\sigma-},Z_{\sigma-})[X^{c},Z^{c}]_{\sigma}^{\eta}+
B_{zz}^\ell(X_{\sigma-},Z_{\sigma-})[Z^{c},Z^{c}]_{\sigma}^{\eta}\\
& \geq \int_{\R^2}\bigg[\frac{b^2}{B_{zz}(X_{\sigma-}+u/\ell,Z_{\sigma-}+v/\ell)}\cdot [Y^c,Y^c]_{\sigma}^\eta\\
&\quad \qquad +B_{zz}(X_{\sigma-}+u/\ell,Z_{\sigma-}+v/\ell)\big[Z^c+(\xi^{(u,v)})^c,Z^c+(\xi^{(u,v)})^c\big]_{\sigma}^\eta\bigg] g(u,v)\mbox{d}u\mbox{d}v\\
&\geq \int_{\R^2} 2b [Y^c,Z^c+(\xi^{(u,v)})^c]_{\sigma}^\eta g(u,v)\mbox{d}u\mbox{d}v\\
&=\int_{\R^2} 2b [Y^c,Z^c]_{\sigma}^\eta g(u,v)\mbox{d}u\mbox{d}v=2b[Y^c,Z^c]_\sigma^\eta=2b[Y,Z]_{\sigma}^\eta.
\end{split}
\end{equation*}
In the first equality above we have used the orthogonality of $X^c$ and $Y^c$ (which gives $d[Y^c,(\xi^{(u,v)})^c]=0$), while the last passage follows from the continuity of paths of $Y$. Therefore, approximating $I_3$ by Riemann sums, we get
$$ I_3\geq 2b[Y,Z]_0^{\tau_n\wedge t}.$$
Putting all the above facts into \eqref{ito} (and noting that the term $I_3$ is multiplied by $1/2$), we see that we have proved the estimate
\begin{equation}\label{interm}
\begin{split}
 \int_{\R^2} \E \Big(\alpha|X_{\tau_n\wedge t}+u/\ell|^p+&\beta|Z_{\tau_n\wedge t}+v/\ell|^{p'}\Big)g(u,v)\mbox{d}u\mbox{d}v \\
&\qquad \qquad \qquad \geq \E B^\ell (X_0,Z_0)+b\E [Y,Z]_0^{\tau_n\wedge t}.
\end{split}
 \end{equation}

\smallskip

\emph{Step 4. Limiting procedure.} Now we let $\ell \to \infty$. Using Lebegue's dominated convergence the theorem, combined with the assumed boundedness of $X$ and $Z$ in $L^p$ and $L^{p'}$, respectively, we get
\begin{align*}
 \int_{\R^2} \E \Big(\alpha|X_{\tau_n\wedge t}+u/\ell|^p+\beta|Z_{\tau_n\wedge t}+v/\ell|^{p'}\Big)g(u,v)\mbox{d}u\mbox{d}v\to  \alpha \|X_{\tau_n\wedge t}\|_p^p+\beta \|Z_{\tau_n\wedge t}\|_{p'}^{p'}.
 \end{align*}
 Furthermore, since $B$ is continuous, we have the almost sure convergence $B^\ell(X_0,Z_0) \to B(X_0,Z_0)$. Hence,  by Fatou's lemma, we get  $\E B(X_0,Z_0)\leq \liminf_{\ell \to \infty} \E B^\ell(X_0,Z_0)$ (the functions $B^\ell$ are nonnegative, so the use of Fatou's lemma is permitted). Plugging these observations into \eqref{interm}, we obtain
 \begin{align*}
 \alpha \|X_{\tau_n\wedge t}\|_p^p+\beta \|Z_{\tau_n\wedge t}\|_{p'}^{p'}&\geq \E B(X_0,Z_0)+b\E [Y,Z]_{0}^{\tau_n\wedge t}\\
&\geq b\mathbb{E}|X_0Z_0|+b\E [Y,Z]_{0}^{\tau_n\wedge t}\\
&\geq b\mathbb{E}Y_0Z_0 +b\E [Y,Z]_{0}^{\tau_n\wedge t}=b\E Y_{\tau_n\wedge t}Z_{\tau_n\wedge t},
 \end{align*}
 where in the second passage we have exploited \eqref{initial} and in the third we used the bound $|Y_0|\leq |X_0|$, which is due to Lemma \ref{l1W}. 
 
 \smallskip
 
 \emph{Step 5. Proof of \eqref{Lp}.} Let us fix an auxiliary positive parameter $\lambda$ and apply the latter inequality to the pair $(\lambda X,\lambda Y)$: this is permitted, since the orthogonality and differential subordination are preserved after such scaling. Dividing both sides by $\lambda$ and optimizing the right-hand side with respect to $\lambda$, we obtain an estimate equivalent to
\begin{align*}
 \E Y_{\tau_n\wedge t}Z_{\tau_n\wedge t}&\leq \frac{(p\alpha)^{1/p}(p'\beta)^{1/p'}}{b} \|X_{\tau_n\wedge t}\|_p\|Z_{\tau_n\wedge t}\|_{p'} \\
&\leq \frac{(p\alpha)^{1/p}(p'\beta)^{1/p'}}{b} \|X\|_p\|Z_{\tau_n\wedge t}\|_{p'}.
\end{align*}
Since $Z$ was taken arbitrarily, this gives
 $$  \| Y_{\tau_n\wedge t}\|_p\leq \frac{(p\alpha)^{1/p}(p'\beta)^{1/p'}}{b} \|X\|_p$$
 and letting $n\to \infty$ yields
 $$  \| Y_{t}\|_p\leq \frac{(p\alpha)^{1/p}(p'\beta)^{1/p'}}{b} \|X\|_p,$$
 by Fatou's lemma. But $t$ was arbitrary, so the claim follows.
 
 \smallskip
 
 \emph{Step 6. Proof of \eqref{Lp'}.} This essentially follows from the symmetry of the conditions \eqref{initial}, \eqref{majorization} and \eqref{det2}. Namely, observe that the function
 $$ \tilde{B}(x,z)=B(z,x)$$
 satisfies 
$$
\tilde{B}(x,z)\geq b|xz| \qquad \mbox{for all }(x,z)\in \R^2,
$$
$$
\tilde{B}(x,z)\leq \beta |x|^{p'}+\alpha|z|^{p}\qquad \mbox{for all }(x,z)\in \R^2
$$
and
$$
\operatorname*{det}D^2\tilde{B}(x,z)\geq b^2 \qquad \mbox{for all }(x,z)\in \bigcup_i D_i.
$$
Therefore, applying \eqref{Lp} with the exponent $p'$, we obtain the estimate
$$ \|Y\|_{p'}\leq \frac{(p'\beta)^{1/p'}(p\alpha)^{1/p}}{b}\|X\|_{p'},$$
which is \eqref{Lp'}.
\end{proof}

An informal, yet important comment is in order. For any $t>0$, let $ [[X,Z]]_t=\operatorname*{limsup} \sum_{j=1}^N |X_{t_j}-X_{t_{j-1}}||Z_{t_j}-Z_{t_{j-1}}|$, where the limit is taken over all $n$ and all partitions $0=t_0<t_1<t_2<\ldots<t_n=t$ of $[0,t]$ with diameter tending to $0$.

\begin{remark}
One can study the estimate \eqref{moment} in the continuous-time context. As proved (directly) by Ba\~nuelos and Wang in \cite{BW,W}, if $Y$ is differentially subordinate to $X$ (the orthogonality is not assumed), then we have the sharp bound
$$ \|Y\|_p\leq (p^*-1)\|X\|_p,\qquad 1<p<\infty.$$
This estimate can be established with a dual argument (cf. \cite{BGO}). For simplicity, let us assume that $X$, $Y$ start from $0$. Then for an arbitrary local martingale $Z$ and a suitable special function $\mathbb{B}$ one can write
$$ \frac{(p^*-1)^p\|X_t\|_p^p}{p}+\frac{\|Z_t\|_{p'}^{p'}}{p'}\geq \E \mathbb{B}(X_t,Z_t)\geq \E [[X,Z]]_t\geq \E [Y,Z]_t=\E Y_tZ_t.$$
 Carrying out the homogenization trick as above, this yields $\|Y_t\|_p\leq (p^*-1)\|X_t\|_p$ and it remains to let $t\to \infty$. The above chain of inequalities is precisely the continuous-time extension of \eqref{chain1} and \eqref{chain2}.
 
It is now instructive to compare this argument with the proof of Theorem \ref{main_tool} we have just presented above. For simplicity, let us assume that $X$, $Y$ have continuous paths and start from zero, and $B$ is of class $C^2$. Furthermore, by homogeneity, we may take $b=1$. Skipping the mollification and the passage to continuous parts of processes, we see that our reasoning yields
\begin{align*}
 \alpha \|X_t\|_p^p+\beta \|Z_t\|_{p'}^{p'}&\geq \E B(X_t,Z_t)\\
&\geq \E [[X,Z+\xi]]_t\geq \E [Y,Z+\xi]_t=\E [Y,Z]_t= \E Y_tZ_t,
\end{align*}
where $\xi_t=\int_0^t a(X_s,Z_s)\mbox{d}X_s$. Thus, in the orthogonal setting, the difference lies  in the possibility of inserting the auxiliary stochastic integral $\xi$ (``an orthogonal compensator''), which is annihilated by the process $Y$ at the end. On the level of convexity of Bellman functions, this corresponds to the possibility of introducing the auxiliary term $a(x,z)$ in \eqref{conv} (compare this to \eqref{conv0}).
\end{remark}

\section{Proof of \eqref{mainin}} 

We turn our attention to the explicit construction of Bellman functions as in Proposition \ref{main_tool}. Since the proposition establishes the bound in $L^p$ and $L^{p'}$ simultaneously (see Remark \ref{dual_remark}), it is enough to consider the case $p>2$.
 
In order to introduce the suitable special function $B$, we distinguish the angular domain
$$ \mathcal{A}=\Big\{(x,z)\in \R^2: |z|>|x|^{p-1}\cdot \mbox{$\sin^{2-p}\frac{\pi}{2p}$}\Big\}$$
and introduce the following auxiliary parameter.

\begin{lemma}
For any $(x,z)\in \mathcal{A}$, there is a unique $\theta=\theta(x,z)\in (\frac{\pi}{2}-\frac{\pi}{2p},\frac{\pi}{2}]$ satisfying the equation
\begin{equation}\label{deftheta}
 |z| \cos^{p-1}\theta=|x|^{p-1}\cos((p-1)(\mbox{$\frac{\pi}{2}$}-\theta)).
\end{equation}
Furthermore, $\theta$ is continuous on $\mathcal{A}$ and of class $C^\infty$ on $\mathcal{A}\setminus \{x=0\}$.
\end{lemma}
\begin{proof}
If $x=0$, then $z\neq 0$ and hence $\theta$ must be the root of the cosine function. This gives $\theta=\pi/2$, as the angle must lie in $(\frac{\pi}{2}-\frac{\pi}{2p},\frac{\pi}{2}]$.  
So, from now on, let us assume that $x\neq 0$. The cosines in \eqref{deftheta} cannot both vanish, so the equation \eqref{deftheta} can be rewritten in the form
\begin{equation}\label{deftheta2}
 \frac{|z|}{|x|^{p-1}}=\frac{\cos((p-1)(\frac{\pi}{2}-\theta))}{\cos^{p-1}\theta}.
\end{equation}
Denote the right-hand side by $F(\theta)$. A direct differentiation shows that
$$ F'(\theta)=\frac{(p-1)\cos((p-2)(
\frac{\pi}{2}-\theta))}{\cos^p\theta},$$
which is positive for $\theta\in [\frac{\pi}{2}-\frac{\pi}{2p},\frac{\pi}{2}]$: indeed,
$$ 0\leq (p-2)(\mbox{$
\frac{\pi}{2}$}-\theta)\leq \frac{(p-2)\pi}{2p}<\frac{\pi}{2}.$$
However, note that $F(\theta)\to \sin^{2-p}\frac{\pi}{2p}$ as $\theta\to \pi/2-\pi/(2p)$ and $F(\theta)\to \infty$ as $\theta\to \pi/2$. Thus we may write
$$ \theta=F^{-1}\left(\frac{|z|}{|x|^{p-1}}\right),$$
with the convention $F^{-1}(\infty)=\pi/2$. Now the regularity of $\theta$ follows from the continuity of $F^{-1}$ on $[\sin^{2-p}\frac{\pi}{2p},\infty]$ and its smoothness on $(\sin^{2-p}\frac{\pi}{2p},\infty)$, direct consequences of the inverse function theorem.
\end{proof}

Consider the parameters
$$ \alpha=\frac{1}{p-1}\cot\frac{\pi}{2p}\sin^{2-p}\frac{\pi}{2p},\qquad \beta=\cot\frac{\pi}{2p}\sin^{(p-2)/(p-1)}\frac{\pi}{2p},\qquad b=\frac{p}{p-1}.$$
We will now prove that the special function
$$ B(x,z)=\begin{cases}
|z|^{p'} & \mbox{if }|z|> |x|^{p-1}\cdot \sin^{2-p}\frac{\pi}{2p}=0,\\
|xz|\big[\tan \theta+\frac{1}{p-1}\tan ((p-1)(\frac{\pi}{2}-\theta))\big] & \mbox{if }|z|>|x|^{p-1}\cdot \sin^{2-p}\frac{\pi}{2p}>0,\\
\alpha |x|^p+ \beta |z|^{p'} & \mbox{if }|z|\leq |x|^{p-1}\cdot \sin^{2-p}\frac{\pi}{2p}
\end{cases}$$
satisfies the requirements listed in Proposition \ref{main_tool}. This will immediately establish the $L^p$ estimate, since
$$ \frac{(p\alpha)^{1/p}(p'\beta)^{1/p'}}{b},$$
the $L^p$ constant given in the proposition, is precisely $\cot \frac{\pi}{2p}$.

\smallskip

The verification of the conditions is rather technical. Let us split the analysis into several parts.

\begin{proof}[Regularity of $B$] Note that $B$ is symmetric with respect to the variables $x$ and $z$, i.e., we have the identity $B(x,z)=B(-x,z)=B(x,-z)$ for all $x$, $z\in \R$. Therefore, it is enough to show that

\medskip

(a) $B$ is continuous on $[0,\infty)^2$;

\smallskip

(b) The partial derivatives $B_x$, $B_z$ exist on $(0,\infty)^2$ and extend to continuous functions on $[0,\infty)^2$, with $B_x(0,z_0)=B_z(x_0,0)=0$ for all $x_0,\,z_0> 0$.

\smallskip

(c) $B$ is of class $C^2$ on the domains $\big\{(x,z)\in (0,\infty)^2: z>x^{p-1}\sin^{2-p}\frac{\pi}{2p}\big\}$ and $\big\{(x,z)\in (0,\infty)^2: z<x^{p-1}\sin^{2-p}\frac{\pi}{2p}\big\}$.
\medskip

 We start with the continuity. It is enough to show that
\begin{equation}\label{lim1}
 \lim_{(x,z)\to (0,z_0)}B(x,z)=B(0,z_0)
\end{equation}
\begin{equation}\label{lim2}
\lim_{(x,z)\to (x_0,y_0)\in \partial\mathcal{A}} B(x,z)=B(x_0,z_0),
\end{equation}
where in both limits the points $(x,z)$ are assumed to lie in $\mathcal{A}\cap [0,\infty)^2$. 
 By \eqref{deftheta}, the formula for $B$ on the set $\{(x,z):|z|>|x|^{p-1}\cdot \sin^{2-p}\frac{\pi}{2p}>0\}$ can be rewritten in the form
$$ B(x,z)=|z|^{p'}\cdot \frac{\sin\theta}{\cos^{p'-1}((p-1)(\frac{\pi}{2}-\theta))}+|xz|\cdot \frac{\tan((p-1)(\frac{\pi}{2}-\theta))}{p-1}.$$
Hence, if $(x,z)\to (0,z_0)$ and $z_0\neq 0$, then $\theta \to \pi/2$ and $B(x,z)\to |z_0|^{p'}$. If $(x,z)\to (0,0)$, then $\theta$ need not stabilize - all we can say that it belongs to $(\frac{\pi}{2}-\frac{\pi}{2p},\frac{\pi}{2}]$ - but then the trigonometric functions above remain bounded, and thus $B(x,z)\to 0=B(0,0)$. This establishes \eqref{lim1}. To show \eqref{lim2}, we may assume that $(x_0,z_0)\neq (0,0)$ (this case is already covered by \eqref{lim1}). Then $\theta\to \frac{\pi}{2}-\frac{\pi}{2p}$ and
$$ B(x,z)\to |x_0z_0|\cdot \frac{p}{p-1}\cot\frac{\pi}{2p}=\alpha |x_0|^p+ \beta |z_0|^{p'},$$
by the very definition of $B$. 

Now we turn our attention to (b): we will present all the details only for $B_x$, the partial derivative with respect to $z$ can be treated similarly. Pick $(x,z)\in \mathcal{A}\cap (0,\infty)^2$. A direct differentiation of \eqref{deftheta2} with respect to $x$ gives
$$ -\frac{(p-1)z}{x^p}=\frac{(p-1)\cos((p-2)(\frac{\pi}{2}-\theta))}{\cos^p\theta}\cdot \theta_x,$$
which combined with \eqref{deftheta2} again, yields
$$ x\theta_x=-\frac{\cos((p-1)(\mb))\cos\theta}{\cos((p-2)(\mb))}.$$
Now, we compute that
\begin{align*} 
&B_x(x,z)\\
&=z\left[\tan\theta+x\theta_x \cos^{-2}\theta\right]+\frac{z}{p-1}\tan((p-1)(\mb ))-zx\theta_x \cos^{-2}((p-1)(\mb ))\\
&=z\tan((p-2)(\mb ))+\frac{z}{p-1}\tan((p-1)(\mb ))-zx\theta_x \cos^{-2}((p-1)(\mb )),
\end{align*}
where the last passage follows from the preceding formula for $x\theta_x$.  
Therefore, if we let $(x,z)\to (0,z_0)$ with $z_0\neq 0$, then $\theta\to \pi/2$ and hence $B_x(x,z)\to 0$. The same is true if  $(x,z)\to (0,0)$, since the trigonometric coefficients in the last formula for $B_x(x,z)$ are bounded: we have $0\leq(p-2)(\pie)\leq (p-1)(\pie)\leq (p-1)\pi/(2p)<\pi/2$. Furthermore, if we take $x,\,z>0$ with $(x,z)\notin \mathcal{A}$ and let $(x,z)\to (0,0)$, then we also obtain $B_x(x,z)=p\alpha x^{p-1} \to 0$. 

Finally, note that for any $(x_0,z_0)\in \partial \mathcal{A}\cap (0,\infty)^2$, we have
\begin{align*}
 \lim_{\mathcal{A}\ni (x,z)\to (x_0,z_0)}B_x(x,z)&=z_0\tan \mbox{$\frac{(p-2)\pi}{2p}$}+\frac{z_0}{p-1}\tan\mbox{$\frac{(p-1)\pi}{2p}$}+z_0\cos^{-1}\mbox{$\frac{(p-2)\pi}{2p}$}\\
&=\frac{p}{p-1}x_0^{p-1}\cdot \sin^{2-p}\cot\mbox{$\frac{\pi}{2p}$}\\
&=\alpha p x_0^{p-1}=\lim_{\mathcal{A}\not \ni (x,z)\to (x_0,z_0)}B_x(x,z).
\end{align*}
It remains to note that the property (c) is trivial (it follows from the smoothness of $\theta$ on $\mathcal{A}\setminus \{x=0\}$, established in the previous lemma).
\end{proof}

\begin{proof}[The conditions \eqref{initial} and \eqref{majorization}]
Again, we may restrict ourselves to $(x,z)\in (0,\infty)^2$.  
Observe that the function
$$ \theta\mapsto \tan\theta+\frac{1}{p-1}\tan((p-1)(\mb ))$$
is increasing on $(\pie,\frac{\pi}{2})$: indeed, its derivative $ \theta\mapsto \cos^{-2}\theta-\cos^{-2}((p-1)(\mb ))$ is positive. In particular, this yields 
$$ \tan\theta+\frac{1}{p-1}\tan((p-1)(\mb ))\geq \tan(\mbox{$\frac{\pi}{2}-\frac{\pi}{2p}$})+\frac{1}{p-1}\tan(\mbox{$\frac{(p-1)\pi}{2p}$})=\frac{p}{p-1}\cot\frac{\pi}{2p},$$
which implies that for $(x,z)\in \mathcal{A}$,
$$ B(x,z)\geq xz\cdot \frac{p}{p-1}\cot\frac{\pi}{2p}\geq bxz.$$
On the other hand, if $(x,z)\notin \mathcal{A}$, then by the Young inequality,
\begin{align*} 
bxz &=\frac{p}{p-1}\cdot \left(x\sin^{(2-p)/p}\mbox{$\frac{\pi}{2p}$}\right)\cdot \left(z \sin^{(p-2)/p}\mbox{$\frac{\pi}{2p}$} \right)\\
&\leq \frac{p}{p-1}\cot\frac{\pi}{2p}\cdot \left(\frac{x^p\sin^{2-p}\frac{\pi}{2p}}{p}+\frac{z^{p'}\sin^{(p-2)/(p-1)}\frac{\pi}{2p}}{p'}\right)=\alpha x^p+\beta z^{p'}=B(x,z)
\end{align*}
and \eqref{initial} is established. To show the majorization, we may assume that $(x,z)\in \mathcal{A}$ (otherwise the claim is trivial: we obtain an equality). Dividing both sides by $x^p$ and using \eqref{deftheta2}, we obtain the equivalent formulation
$$ L(\theta):=F(\theta) \left[\tan\theta+(p-1)^{-1}\tan((p-1)(\mb ))\right]-\alpha -\beta (F(\theta))^{p'}\leq 0.$$
Both sides are equal for $\theta=\frac{\pi}{2}-\frac{\pi}{2p}$ (this corresponds to the majorization for $(x,z)\in \partial \mathcal{A}$), so it is enough to show that $L'(\theta)\leq 0$ for $\theta\in (\frac{\pi}{2}-\frac{\pi}{2p},\frac{\pi}{2})$. We compute that
\begin{align*}
 L'(\theta)&=F'(\theta)\left[\tan\theta+(p-1)^{-1}\tan((p-1)(\mb ))-\beta p'(F(\theta))^{p'-1}\right]\\
 &\quad +F(\theta)\left[\cos^{-2}\theta-\cos^{-2}((p-1)(\mb ))\right].
 \end{align*}
Multiplying both sides by $\cos^{p+1}\theta\cos((p-1)(\pie))/\cos((p-2)(\pie))$ and using the identity 
$$\cos^{-2}\theta-\cos^{-2}((p-1)(\mb ))=\cos((p-2)(\pie))\cos(p(\pie)),$$
we rewrite the inequality $L'(\theta)\leq 0$ in the form
\begin{align*}
 &(p-1)\sin\theta\cos((p-1)(\mb ))+\sin((p-1)(\mb ))\cos\theta\\
&\qquad -\beta p\cos^{p'}((p-1)(\mb ))+\cos(p(\mb ))\leq 0,
\end{align*}
which is further equivalent to $ \sin \theta-\beta \cos^{p'-1}((p-1)(\mb ))\leq 0,$ 
or
\begin{equation}\label{equiform}
  \sin\theta \cos^{1-p'}((p-1)(\mb ))\leq \beta.
\end{equation}
Both sides are equal for $\theta=\frac{\pi}{2}-\frac{\pi}{2p}$; denoting the left-hand side by $M(\theta)$, we have
$$ M'(\theta)=-\sin((p-2)(\mb ))\cos^{-p'}((p-1)(\mb ))\leq 0.$$
This shows that $M$ is decreasing on $(\frac{\pi}{2}-\frac{\pi}{2p},\frac{\pi}{2})$, which confirms \eqref{equiform} and establishes the desired majorization.
\end{proof}

\begin{proof}[The convexity of $B$ and the condition \eqref{det2}]
Clearly, we may restrict ourselves to $x,\,z>0$. We consider two possibilities separately.

\smallskip

\emph{The case $z<x^{p-1}\cdot \sin^{2-p}\frac{\pi}{2p}$.} We have
$$ B_{xx}(x,z)=\alpha p(p-1)x^{p-2},\quad B_{xz}(x,z)=0,\quad B_{zz}(x,z)=\beta p'(p'-1)z^{p'-2},$$
which implies
$$ \operatorname*{det}D^2B(x,z)=pp'\alpha \beta x^{p-2}z^{p'-2}>\left(\frac{p}{p-1}\right)^2 \cot^2\frac{\pi}{2p}>b^2.$$

\smallskip

\emph{The case $z>x^{p-1}\cdot \sin^{2-p}\frac{\pi}{2p}$.} To keep the notation shorter, introduce the function
$$ C(\theta)=\tan\theta+\frac{1}{p-1}\tan((p-1)(\mb) ).$$
As we have already seen above, the differentiation of \eqref{deftheta2} with respect to the variable $x$ leads to the identities
\begin{equation}\label{thx}
 x\theta_x=-\frac{\cos((p-1)(\mb))\cos\theta}{\cos((p-2)(\mb))}
\end{equation}
and
\begin{equation}\label{thz}
 z\theta_z=\frac{\cos((p-1)(\mb))\cos\theta}{(p-1)\cos((p-2)(\mb))}.
\end{equation}
Let us also note that
\begin{equation}\label{diffC}
C'(\theta)=\frac{\cos\big((p-2)(\mb)\big)\cos(p(\mb))}{\cos^2\theta\cos^2((p-1)(\mb))}
\end{equation}

Now, since $B(x,z)=xzC(\theta)$, we get $ B_x(x,z)=zC(\theta)+xzC'(\theta)\theta_x,$ 
which, after some calculations involving \eqref{thx} and \eqref{diffC}, takes the much nicer form
\begin{equation}\label{bex}
 B_x(x,z)=\frac{pz}{p-1}\tan((p-1)(\mb)).
\end{equation}
Similarly, we check that
\begin{equation}\label{bez}
 B_z(x,z)=xC(\theta)+xzC'(\theta)\theta_z=\frac{px}{p-1}\tan\theta.
\end{equation}
Consequently, differentiating again and using $x\theta_x=-(p-1)z\theta_z$ (which follows at once from \eqref{thx} and \eqref{thz}), we compute that
\begin{equation}\label{bexx}
 B_{xx}(x,z)=-pz\cos^{-2}((p-1)(\mb))\theta_x=\frac{pz}{x}\cdot \frac{(p-1)z\theta_z}{\cos^2((p-1)(\mb))},
\end{equation}
\begin{equation}\label{bezz}
 B_{zz}(x,z)=\frac{px}{(p-1)z}\cdot \frac{z\theta_z}{\cos^2\theta}
\end{equation}
and 
\begin{align*}
 B_{xz}(x,z)&=\frac{p}{p-1}\tan\theta+\frac{px}{p-1}\cdot \frac{\theta_x}{\cos^2\theta}\\
 &=\frac{pz\theta_z}{\cos^2\theta}\left[\frac{\sin \theta\cos((p-2)(\mb))}{\cos((p-1)(\mb))}-1\right]=\frac{pz\theta_z}{\cos\theta}\cdot \frac{\sin((p-2)(\mb))}{\cos((p-1)(\mb))}.
 \end{align*}
 The above formulas allow us to derive that
\begin{align*}
\operatorname*{det}D^2B(x,z)&=\frac{p^2(z\theta_z)^2}{\cos^2\theta\cos^2((p-1)(\mb))}\left(1-\sin^2((p-2)(\mb))\right)\\
&=\frac{p^2(z\theta_z)^2\cos^2((p-2)(\mb))}{\cos^2\theta\cos^2((p-1)(\mb))}=\left(\frac{p}{p-1}\right)^2=b^2,
\end{align*}
 as desired. Finally, let us observe that in both cases above the inequality $B_{xx}>0$ holds; combining this with the fact that $B$ is of class $C^1$ and symmetric with respect to the variables $x$ and $z$, we get the global convexity of $B$.
\end{proof}

\section{The Hilbert and Riesz transforms}  Let $d\geq 1$  and denote by $C^{\infty}_0(\R^d)$ the smooth real-valued functions of compact support defined on the $d$-dimensional Euclidean space $\R^d$. For $j=1, 2, \dots, d$ and $f\in C^{\infty}_0(\R^d)$, the classical Riesz transforms on $\R^d$ are the singular integrals operators defined   by 
\begin{align}\label{RieszT}
    R_{j}f(x) = c_d\,\, p.v.\int_{\R^d}\frac{y_j}{|y|^{d+1}}f(x-y)  \mbox{d}y,\qquad c_d=\pi^{-\tfrac{d+1}{2}}\Gamma(\tfrac{d+1}{2}). 
    \end{align}
 Equivalently, as Fourier multipliers, they are given  by
 \begin{equation}\label{RieszForierM}
\widehat{{R_{j}f}}(\xi)=\frac{i\xi_j}{|\xi|}\widehat{f}(\xi).  
\end{equation}
When $d=1$ both definitions reduce to the Hilbert transform defined  by 
 \begin{align*}
Hf(x)=\frac{1}{\pi}\, p.v.\int_{\R} \frac{f(x-y)}{y} \mbox{d}y
\end{align*} 
and 
\begin{equation}\label{Fouier:MulH}
\widehat{Hf}(\xi)=i\sign(\xi)\widehat{f}(\xi). 
\end{equation}

The sharp $L^p$ constant for the Hilbert transform, $ \cot\left(\frac{\pi}{2p^*}\right)$, was obtained by Pichorides in \cite{Pic}.  In  \cite{IwaMar} it is shown that the classical Calder\'on-Zygmund method or rotations, in combination with de Leeuw's Fourier extension theorem \cite[Theorem 2.5.15]{Gra}, immediately gives the same sharp $L^p$ constant for the Riesz transforms in all dimensions.  Pichorides' argument is purely analytic and relies on the fact that the composition of the (direct) Bellman function with a harmonic function $u$ and its conjugate $v$ 
 is superharmonic.

The dual martingale inequalities \eqref{Lp} and \eqref{Lp'} immediately yield the corresponding dual inequalities for the Riesz transforms via their Gundy-Varopoulos representation, as in \cite{BW}. The aim of this section is to show the dual formulation of the sharp inequalities for the Riesz transforms directly from the dual Bellman function introduced above, thereby providing a purely analytic argument that avoids any appeal to martingale inequalities. The same proof also yields the dual inequality for  the Hilbert transform.

In what follows, we will always work with functions in $C^{\infty}_0(\R^d)$. Classical density arguments extend the inequalities to general functions in $L^p(\R^d)$. 

\begin{theorem}\label{anal:proof}
For every \(f\in C_0^\infty(\R^d)\),
\begin{equation}\label{Lp-anal}
\|R_jf\|_{L^p(\R^d)}
\le
\frac{(p\alpha)^{1/p}(p'\beta)^{1/p'}}{b}\,
\|f\|_{L^p(\R^d)}
\end{equation}
and 
\begin{equation}\label{Lp-anal2}
\|R_jf\|_{L^{p'}(\R^d)}
\le
\frac{(p\alpha)^{1/p}(p'\beta)^{1/p'}}{b}\,
\|f\|_{L^{p'}(\R^d)}.
\end{equation}

\end{theorem} 
These are the same inequalities for the Riesz transforms as the inequalities \eqref{Lp} and \eqref{Lp'} for martingales. As in Remark \ref{dual_remark}, we assume that $p>2$ in what follows.

First, for  $f\in C^{\infty}_0(\R^d)$, let $P_yf(x)$ be the convolution of $f$ with the Poisson kernel 
\[P_y(x)=\frac{c_d\, y}{\left(|x|^2+|y|^2\right)^{\frac{d+1}{2}}}, \quad x\in \R^d, \quad y>0, 
\]
where $c_d$ is the constant as in the  definition of the Riesz transforms. We recall the following well-known Littlewood-Paley identity. 

\begin{lemma}\label{Littlewood-Paley:1} For all  $f, g\in C^{\infty}_0(\R^d)$, set $u(x, y)=P_yf(x)$, $v(x, y)=P_yg(x)$  and fix $j\in \{1, \dots, d\}$. Then 
\begin{equation*}\label{eq:Riesz-area}
\int_{\mathbb R^d} R_j f(x)\,g(x)\,\mbox{d}x
=
2\int_0^\infty y\Big(\int_{\mathbb R^d}
u_{x_j}(x,y)\, v_y(x,y)
-\partial_y u(x,y)\,\partial_{x_j} v(x,y)dx\Big)\,\mbox{d}y.
\end{equation*}
\end{lemma} 
\begin{proof} 
To verify this identity, consider the Fourier transform of $h\in C^{\infty}_0(\mathbb{R}^d)$ and its inverse:
\[
\widehat{h}(\xi)=\int_{\mathbb R^d} h(x)e^{2\pi i x\cdot \xi}\,\mbox{d}x,
\qquad
h(x)=\int_{\mathbb R^d}\widehat h(\xi)e^{-2\pi i x\cdot \xi}\,\mbox{d}\xi.
\]

 Using Placherel's theorem on the integral over $\R^d$ and the fact that 
\[
\widehat{\partial_{x_j} h}(\xi)=-\,2\pi i\,\xi_j\,\widehat h(\xi), \quad \widehat{P_y}(\xi)=e^{-\pi y|\xi|}, \quad 
\widehat{R_j h}(\xi)= \frac{i\,\xi_j}{|\xi|}\,\widehat h(\xi), 
\]
the identity follows after integrating with respect to $y$.  \end{proof} 
Our goal is now to prove the following Proposition 

\begin{proposition}\label{prop:key} For all  $f, g\in C^{\infty}_0(\R^d)$  and $j=1,\,2,\,\ldots,\, d$, we have

\begin{equation*}
2b\int_0^\infty\int_{\R^d}
\bigl|u_{x_j}(x,y)\, v_y(x,y)
-\partial_y u(x,y)\,\partial_{x_j} v(x,y)\bigr|
\,y\,\mbox{d}x\,\mbox{d}y
\le
\int_{\R^d}B(f(x),g(x))\,\mbox{d}x. 
\end{equation*}
\end{proposition} 

Fix $f\in C_0^{\infty}(\R^d)$. Assuming the Proposition we have, using the majorization \eqref{majorization}  of the  Bellman function $B$ and Lemma \ref{Littlewood-Paley:1}, 
\[
b\left|\int_{\R^d}R_jf(x)\,g(x)\,\mbox{d}x\right|
\le
\alpha\int_{\R^d}|f(x)|^p\,\mbox{d}x+\beta\int_{\R^d}|g(x)|^{p'}\,\mbox{d}x, 
\]
for every \(g\in C_0^\infty(\R^d)\). Replacing \(f\) by \(\lambda f\) and \(g\) by \(\lambda g\), with $\lambda>0$,
we have 
\[
b\lambda^2\left|\int_{\R^d}R_jf(x)\,g(x)\,\mbox{d}x\right|
\le
\alpha \lambda^p\|f\|_{L^p}^p+\beta \lambda^{p'}\|g\|_{L^{p'}}^{p'}.
\]
Optimizing in $\lambda$  yields
\[\left|\int_{\mathbb R^d} R_j f(x)\,g(x)\,\mbox{d}x\right|
\le
\frac{(p\alpha)^{1/p}(p'\beta)^{1/p'}}{b}\,
\|f\|_{L^p}\|g\|_{L^{p'}}.
\]
Taking the supremum over all \(g\) with \(\|g\|_{L^{p'}}\le1\) completes the proof of Theorem \ref{anal:proof}. 

We now turn to the proof of Proposition \ref{prop:key} which is presented in several lemmas.  
Consider the mollified Bellman  function \(B^\ell\) as in \eqref{Mollified:1}. Since $B$ satisfies the requirements listed in Proposition \ref{main_tool}, in particular $B$ is continuous on $\R^2$, it follows from the standard convolution arguments that 
\begin{equation}\label{limB:ell}
B^\ell(x,z)\to B(x,z)
\qquad\text{for every }(x,z)\in\R^2. 
\end{equation} 
Furthermore,  if \(K\subset \R^2\) is compact,  the convergence is uniform on \(K\).
\begin{lemma}\label{prop:mollified-B}
For every \(\ell\), \(B^\ell\) satisfies
\begin{equation}\label{eq:mollified-hessian}
D^2B^\ell(x, z)\ge0,\quad \text{and}
\quad
\det D^2B^\ell(x, z)\ge b^2, 
\quad(x, z)\in \R^2.
\end{equation}

\end{lemma}
\begin{proof} It follows from \eqref{i1}, \eqref{i2}, \eqref{i3} that
\begin{equation*}
D^2B^\ell(x,z)=\int_{\R^2} D^2B(x+u/\ell,z+v/\ell)\,g(u,v)\,du\,dv,
\end{equation*}
or more explicitly, in matrix notation, $D^2B^\ell(x,z)$ equals
\begin{align*}
\begin{pmatrix}
\displaystyle \int_{\bR^2}B_{xx}(x+u/\ell,z+v/\ell)g(u,v)\,\mbox{d}u\,\mbox{d}v
&
\displaystyle \int_{\bR^2} B_{xz}(x+u/\ell,z+v/\ell)g(u,v)\,\mbox{d}u\,\mbox{d}v
\\[1.2em]
\displaystyle \int_{\bR^2} B_{xz}(x+u/\ell,z+v/\ell)g(u,v)\,\mbox{d}u\,\mbox{d}v
&
\displaystyle \int_{\bR^2} B_{zz}(x+u/\ell,z+v/\ell)g(u,v)\,\mbox{d}u\,\mbox{d}v
\end{pmatrix}.
\end{align*}
Therefore
\[
\det D^2B^\ell(x,z)
=
\left(\int B_{xx}g\right)\left(\int B_{zz}g\right)
-
\left(\int B_{xz}g\right)^2,
\]
where, for simplicity, the entries \(B_{xx},B_{xz},B_{zz}\) are evaluated at
\((x+u/\ell,z+v/\ell)\) and all integrals are over $\R^2$. 
Let us consider  the function 
$$F(A)=(\det(A))^{1/2}$$ 
in the space (cone) of positive semidefinite $2\times 2$ matrices. 
From \cite[Corollary II.3.21, p.47]{Bha}, $F$ is concave  and it follows by Jensen's inequality that 

\[
\left[
\left(\int B_{xx}g\right)\left(\int B_{zz}g\right)
-
\left(\int B_{xz}g\right)^2
\right]^{1/2}
\ge
\int
\left(B_{xx}B_{zz}-B_{xz}^2\right)^{1/2} g.
\]
Thus, since
\[
B_{xx}B_{zz}-B_{xz}^2=\det D^2B\ge b^2
\]
almost everywhere, we get
\[
\left(\det D^2B^\ell(x,z)\right)^{1/2}
\ge
\int b\,g(u,v)\,\mbox{d}u\,\mbox{d}v=b.
\]
Consequently, we have \(
\det D^2B^\ell(x,z)\ge b^2,
\)
 as desired.
\end{proof} 

We next recall the following elementary lemma from linear algebra and give its proof for completeness. 
\begin{lemma}\label{lem:trace-det-2x2}
Let $A$ and $S$ be $2\times 2$ positive semidefinite matrices. Then
\begin{equation}\label{eq:trace-det}
\operatorname{tr}(AS)\ge 2\sqrt{\det A\,\det S}.
\end{equation}
\end{lemma}

\begin{proof}
Let
\[
B=A^{1/2}SA^{1/2}.
\]
Then $B$ is positive semidefinite and
$
\operatorname{tr}(B)=\operatorname{tr}(AS).
$
Thus, 
\[
\det(B)=\det(A^{1/2})\,\det(S)\,\det(A^{1/2})
=\det(A)\det(S).
\]
Let $\lambda_1,\lambda_2\ge 0$ be the eigenvalues of $B$. Then
\[
\operatorname{tr}(AS)=\operatorname{tr}(B)=\lambda_1+\lambda_2
\ge 2\sqrt{\lambda_1\lambda_2}
=2\sqrt{\det(B)}
=2\sqrt{\det A\,\det S},
\]
which verifies the inequality.  
\end{proof}

Returning to the notation of Lemma \ref{Littlewood-Paley:1}, 
for a fixed  \(j\in\{1,\,2,\,\ldots,\,d\}\), define
\[
\xi_j=\binom{\partial_{x_j}u}{\partial_{x_j}v},
\qquad
\eta=\binom{\partial_yu}{\partial_yv},
\qquad
S_j:=\xi_j\xi_j^T+\eta\eta^T.
\]
That is, 
\[
S_j=
\begin{pmatrix}
u_{x_j}^2+u_y^2 & u_{x_j}v_{x_j}+u_yv_y\\[0.4em]
u_{x_j}v_{x_j}+u_yv_y & v_{x_j}^2+v_y^2
\end{pmatrix}.
\]
Then, 
\begin{equation}\label{detS_j}
\det S_j=
\bigl(u_{x_j}v_y-u_yv_{x_j}\bigr)^2.
\end{equation} 
Now, if \(A=(a_{mn})\) is a \(2\times2\) matrix and \(w\in\R^2\), then
\[
\langle Aw,w\rangle=\tr(Aww^T).
\]
Thus applying this first with \(w=\xi_j\), and then with \(w=\eta\), and adding, gives
\begin{align}
\langle A\xi_j,\xi_j\rangle+\langle A\eta,\eta\rangle
&=\tr(A\xi_j\xi_j^T)+\tr(A\eta\eta^T)\\
&=\tr(A(\xi_j\xi_j^T+\eta\eta^T))
=\tr(AS_j).\nonumber
\end{align}
Define 
\[
\Phi^\ell(x,y):=B^\ell(u(x,y),v(x,y)).
\]
Since \(B^\ell\in C^\infty(\R^2)\), the function \(\Phi^\ell\) belongs to  $C^\infty(\Rp)$.

\begin{lemma}\label{lem:chain}
Let
\[
A^\ell=D^2B^\ell(u(x,y),v(x,y)).
\]
Then
\[
\Delta \Phi^\ell(x,y)
=
\sum_{k=1}^d\langle A^\ell\xi_k,\xi_k\rangle+\langle A^\ell\eta,\eta\rangle.
\]
In particular, for each fixed \(j\),
\begin{equation}\label{lboundtrace}
\Delta\Phi^\ell(x,y)\ge \langle A^\ell\xi_j,\xi_j\rangle+\langle A^\ell\eta,\eta\rangle
=\tr(A^\ell S_j).
\end{equation}
\end{lemma}
\begin{proof}
For \(k\in\{1,\,2,\,\ldots,\,d\}\), we have
\[
(\Phi^\ell)_{x_k}
=
(B^\ell)_x(u,v)u_{x_k}+(B^\ell)_z(u,v)v_{x_k}.
\]
Differentiating once more,
\[
(\Phi^\ell)_{x_kx_k}
=
\langle A^\ell\xi_k,\xi_k\rangle
+(B^\ell)_x(u,v)u_{x_kx_k}
+(B^\ell)_z(u,v)v_{x_kx_k}.
\]
Similarly,
\[
(\Phi^\ell)_y
=
(B^\ell)_x(u,v)u_y+(B^\ell)_z(u,v)v_y,
\]
and hence
\[
(\Phi^\ell)_{yy}
=
\langle A^\ell\eta,\eta\rangle
+(B^\ell)_x(u,v)u_{yy}
+(B^\ell)_z(u,v)v_{yy}.
\]
Summing over \(k=1,\,2,\,\ldots,\,d\) and adding the \(y\)-term, we obtain
\[
\Delta\Phi^\ell
=
\sum_{k=1}^d\langle A^\ell\xi_k,\xi_k\rangle+\langle A^\ell\eta,\eta\rangle
+(B^\ell)_x(u,v)\Delta u+(B^\ell)_z(u,v)\Delta v.
\]
Since \(u\) and \(v\) are harmonic in \(\Rp\), \(\Delta u=\Delta v=0\), and the first formula follows.
The second follows by discarding the nonnegative terms \(\langle A^\ell\xi_k,\xi_k\rangle\) for \(k\neq j\),
since \(A^\ell\ge0\) by Lemma~\ref{prop:mollified-B}.
\end{proof}

Since \(A^\ell\ge0\) and \(S_j\ge0\), the Lemma together with \eqref{eq:trace-det} and \eqref{detS_j} immediately give
\begin{equation}\label{eq:pointwise}
\Delta\Phi^\ell(x,y)
\ge 2\sqrt{\det A^\ell\,\det S_j}\geq 
2b\,\bigl|u_{x_j}(x,y)v_y(x,y)-u_y(x,y)v_{x_j}(x,y)\bigr|.
\end{equation}
From this and Lemma \ref{Littlewood-Paley:1} we have 

\begin{corollary}\label{3:4;there} For all  $f, g\in C^{\infty}_0(\R^d)$, 
\begin{align*}
\bigl|\int_{\mathbb R^d} R_j f(x)\,g(x)\,\mbox{d}x\bigl|
& \leq 
2\int_0^\infty \int_{\mathbb R^d}
\big|u_{x_j}(x,y)\, v_y(x,y)
-\partial_y u(x,y)\,\partial_{x_j} v(x,y)\bigl|y\, \mbox{d}x\,\mbox{d}y\\
&\leq \frac{1}{b}\int_0^\infty \int_{\mathbb R^d}\Delta\Phi^\ell(x,y) y\, \mbox{d}x\,\mbox{d}y.
\end{align*}
\end{corollary}

In view of \eqref{limB:ell}, the proof of Proposition \ref{prop:key} will be complete once we show that

\begin{equation}\label{Ste;identity}
\frac{1}{b}\int_0^\infty \int_{\mathbb R^d}\Delta\Phi^\ell(x,y) y\, dx\,dy=\frac{1}{b} \int_{\mathbb R^d}\Phi^\ell(x,0)dx,
\end{equation} 
since the quantity on the right equals 
$$\frac{1}{b}\int_{\mathbb R^d} B^{\ell}(f(x), g(x))\mbox{d}x. 
$$ 

The equality \eqref{Ste;identity} would follow from \cite[Lemma 2, p.~87]{Ste-SI}, provided that the function $\Phi^\ell(x,y)$ satisfies the necessary decay assumptions of the lemma. More precisely, we require that $\Phi^\ell \in C^2(\mathbb{R}^{d+1}_+)\cap C(\overline{\mathbb{R}^{d+1}_+})$ and decays suitably at infinity. This will be achieved by a small modification of the function $\Phi^\ell$, as described in the following lemma.

\begin{lemma}
Let \[
\widetilde\Phi^\ell(x,y)
=
\widetilde B^\ell(u(x,y),v(x,y)), 
\]
where 
\begin{align*}
&\widetilde B^\ell(u(x, y),\,v(x, y))\\
&=
B^\ell(u(x, y), v(x, y))-B^\ell(0,0)-B_x^\ell(0,0)\,u(x, y)-B_z^\ell(0,0)\,v(x, y). 
\end{align*}
Then $\widetilde\Phi^\ell\in C^2(\R^{d+1}_+)\cap C(\overline{\R^{d+1}_+})$,
$\widetilde\Phi^\ell$ is suitably small at infinity in the sense of
Stein's Lemma, and
\begin{equation}\label{eq:stein-identity}
\int_0^\infty\int_{\R^d}
\Delta \Phi^\ell(x,y)\,y\,\mbox{d}x\,\mbox{d}y
=
\int_{\R^d} B^\ell(f(x),g(x))\,\mbox{d}x
+ E_\ell,
\end{equation}
where the error $E_\ell\to 0$ as $\ell\to\infty$.
\end{lemma}

\begin{proof}
Since $B^\ell\in C^\infty(\R^2)$ and $u,v$ are harmonic in
$\R^{d+1}_+$, both $\Phi^\ell$ and $\widetilde\Phi^\ell$ are $C^\infty$
in $\R^{d+1}_+$. Moreover, since $u(\cdot,y)\to f$ and $v(\cdot,y)\to g$
uniformly on compact sets as $y\downarrow 0$, both $u$ and $v$ extend continuously to
$\overline{\R^{d+1}_+}$. In addition, the fact that $u, v$ are harmonic gives 
\[
\Delta \widetilde\Phi^\ell
=
\Delta \Phi^\ell.
\]
By  construction,
\[
\widetilde B^\ell(0,0)=0,\qquad \nabla \widetilde B^\ell(0,0)=0.
\]
Since $\widetilde B^\ell$ is of class $C^2$, the Taylor expansion gives
\[
|\widetilde B^\ell(s,t)|\le C_\ell (s^2+t^2),
\qquad
|\nabla \widetilde B^\ell(s,t)|\le C_\ell (|s|+|t|).
\]
For $f,g\in C_0^\infty(\R^d)$, their Poisson extensions satisfy, for
large $|x|+y$,
\[
|u(x,y)|+|v(x,y)|\le C_N (1+|x|+y)^{-d},
\]
\[
|\nabla u(x,y)|+|\nabla v(x,y)|
\le C_N (1+|x|+y)^{-d-1}.
\]
Therefore
\begin{align}
|\widetilde\Phi^\ell(x,y)|
\le C_\ell (1+|x|+y)^{-2d},
\qquad
|\nabla \widetilde\Phi^\ell(x,y)|
\le C_\ell (1+|x|+y)^{-2d-1}.
\end{align}
With these estimates, 
Stein’s identity \cite[p.~87]{Ste-SI} applies and we have
\[
\int_0^\infty\int_{\R^d}
\Delta \widetilde\Phi^\ell(x,y)\,y\,\mbox{d}x\,\mbox{d}y
=
\int_{\R^d} \widetilde\Phi^\ell(x,0)\,\mbox{d}x.
\]
Since $\Delta \widetilde\Phi^\ell=\Delta \Phi^\ell$, this gives
\[
\int_0^\infty\int_{\R^d}
\Delta \Phi^\ell(x,y)\,y\,\mbox{d}x\,\mbox{d}y
=
\int_{\R^d} \widetilde B^\ell(f(x),g(x))\,\mbox{d}x.
\]

Let \(K=\operatorname{supp} f\cup \operatorname{supp} g\).
If \(x\notin K\), then \(f(x)=g(x)=0\), and hence
\[
\widetilde B^\ell(f(x),g(x))=\widetilde B^\ell(0,0)=0.
\]
Therefore the integrand vanishes outside \(K\), and
\[
\int_{\R^d}\widetilde B^\ell(f,g)\,\mbox{d}x
=
\int_K \widetilde B^\ell(f,g)\,\mbox{d}x.
\]
Thus
\begin{align*}
\int_{\R^d}\widetilde B^\ell(f,g)\,\mbox{d}x
&=
\int_K B^\ell(f,g)\,\mbox{d}x
-
B^\ell(0,0)|K|\\
&\quad
-
B_x^\ell(0,0)\int_{\R^d}f\,\mbox{d}x
-
B_z^\ell(0,0)\int_{\R^d}g\,\mbox{d}x.
\end{align*}
Since $B(0,0)=0$ and $B(x,z)=O(|x|^p+|z|^{p'})$ near $(0,0)$, the standard
properties of convolutions  give
\[
B^\ell(0,0)\to 0,\qquad
B_x^\ell(0,0)\to 0,\qquad
B_z^\ell(0,0)\to 0.
\]
Thus the error term
\[
E_\ell
:=
-
B^\ell(0,0)|K|
-
B_x^\ell(0,0)\int_{\R^d}f\,\mbox{d}x
-
B_z^\ell(0,0)\int_{\R^d}g\,\mbox{d}x
\]
satisfies $E_\ell\to 0$, as $\ell\to\infty$.

Combining the above identities yields \eqref{eq:stein-identity}.
This completes the proof of \eqref{Ste;identity} and hence the proof of Proposition \ref{prop:key}. \end{proof}
\bigskip

{\it  AI Disclosure:} This paper was prepared nearly a year ago, and no AI tools were used.
\bigskip
%
%
%

\end{document}